\documentclass[12pt]{article}
\usepackage[dvips]{color}
\usepackage{amssymb}

\newcommand{\g}{\gamma}
\newcommand{\G}{{\mathit{\Gamma}}}

\newcommand{\W}{{\mathit{\Omega}}}
\newcommand{\w}{\omega}

\newcommand{\f}{\varphi}
\renewcommand{\l}{\lambda}
\renewcommand{\L}{{{\mathit\Lambda}}}
\newcommand{\z}{\zeta}

\newcommand{\C}{\mathbb C}
\newcommand{\R}{\mathbb R}
\newcommand{\Q}{\mathbb Q}
\newcommand{\Z}{\mathbb Z}

\newcommand{\Ps}{\mathbb P}
\newcommand{\B}{\mathbb B}
\renewcommand{\H}{\mathbb H}
\renewcommand{\S}{\mathbb S}

\newcommand{\la}{\langle}
\newcommand{\ra}{\rangle}
\newcommand{\tr}{\;^t}
\newcommand{\lto}{\longrightarrow}
\newcommand{\mod}{\ {\rm mod}\ }
\newcommand{\qed}{ {\rm [\kern-.4pt]}  \bigbreak }
\newcommand{\qedd}{ \quad {\rm [\kern-.4pt]}  }
\def \overset#1#2{{\buildrel #1\over#2}}
\def \leftset#1#2{{\scriptstyle #1}\!#2\quad }
\def \comment#1{{ }}
\newtheorem{theorem}{Theorem}[section]
\newtheorem{proposition}{Proposition}[section]

\newtheorem{fact}{Fact}[section]
\title{Theta constants associated with the cyclic triple coverings 
of the complex projective line branching at six points 
\footnote{
$1991$ {\em Mathematics Subject Classification\/} 
Primary 32N15; Secondary 11F55, 14J15, 32G20. 
}}
\author{Keiji Matsumoto}
\date{ July 21, 2000 }
\begin{document}
\maketitle
\begin{abstract}
Let $\psi$ be the period map for a family of the cyclic triple coverings of 
the complex projective line branching at six points.    
The symmetric group $S_6$ acts on this family and on its image under $\psi.$
In this paper, we give an $S_6$-equivariant expression of $\psi^{-1}$ 
in terms of fifteen theta constants. 
\end{abstract}
\section{Introduction}
Let $C(\l)$ be the cyclic triple covering of the complex projective line 
$\Ps^1$ branching at six points $\l_1,\dots,\l_6:$ 
$$C(\l): \ w^3=\prod_{i=1}^6 (z-\l_i).$$ 
The moduli space of such curves with a homology marking 
can be regarded as the configuration space $\L$ of ordered six 
distinct points on $\Ps^1,$ which is defined by 
$$GL_2(\C)\backslash \{\l=(\l_{ij})\in M(2,6)\mid \l\la ij\ra=
\left|\matrix{\l_{1i}&\l_{1j}\cr \l_{2i}&\l_{2j}\cr}\right| \ne 0\}/(\C^*)^6.$$
Note that the symmetric group $S_6$ naturally acts on $\L.$
It is shown in \cite{Yo2} that the map  
$$\iota:\L\ni \l\mapsto 
[\dots,y_{\la ij; kl;mn\ra},\dots]=
[\dots,\l\la ij\ra \l\la kl\ra \l\la mn\ra,\dots]\in \Ps^{14}$$
is an $S_6$-equivariant embedding and that its image 
is an open subset of $Y$ defined by linear and cubic equations.  

The normalized period matrix $\Omega$ of $C(\l)$ 
with a homology marking belongs to 
the Siegel upper half space $\S^4$ of degree $4.$ 
By our assignment of a homology marking, $\W$ can be identified with 
an element of $3$-dimensional complex ball 
$\B^3=\{x\in \Ps^3\mid \tr\bar xHx<0\},$ 
where $H={\rm diag}(1,1,1,-1).$  
In this way, we get a multi-valued map $\psi:\L \to \B^3\subset \S^4,$ which 
is called the period map. 
Results in \cite{DM} and \cite{Ter} imply that 
the image of $\psi$ is an open dense subset of $\B^3,$ 
the monodromy group of $\psi$ is the principal congruence subgroup $\G(1-\w)$
of level $(1-\w)$ of 
$\G=\{g\in GL_4(\Z[\w])\mid\tr\bar gHg=H\},$
and that the inverse of $\psi$ is single valued. 

In this paper, we express the inverse of the period map $\psi$ 
in terms of fifteen theta constants. More precisely, for the two isomorphisms 
$\psi:\L\to \psi(\L)/\G(1-\w)$ and 
$\iota:\L\to \iota(\L)\subset Y\subset \Ps^{14},$ 
we present an isomorphism 
${\mathit\Theta}:\psi(\L)/\G(1-\w)\to \iota(\L)$ such that 
the following diagram commutes:
\begin{equation}
\label{diagram}
\matrix{ \L & \overset\psi\lto & \psi(\L)/\G(1-\w) \cr
            &                  &               \cr
  \quad \leftset\iota\downarrow& \leftset{\mathit\Theta}\swarrow & \cr
            &                  &               \cr  
 \iota(\L)\subset Y\subset \Ps^{14}.  &    & \cr}
\end{equation}
The map ${\mathit\Theta}$ is given by the ratio of the cubes of 
the fifteen theta constants on $\S^4$  
which are invariant under the action of $\G(1-\w)$ embedded in $Sp(8,\Z).$
In particular, linear and cubic relations among the cubes of fifteen theta 
constants coincide with the defining equations of $Y\subset \Ps^{14}.$

It is known that $\G/\la\G(1-\w),-I_4\ra$ is isomorphic to $S_6,$  
which naturally acts on $\psi(\L)/\G(1-\w).$ 
The period map $\psi$ is $S_6$-equivariant. 
By considering the action $S_6\simeq\G/\la\G(1-\w),-I_4\ra$ 
on the fifteen theta characteristics,
we label fifteen theta constants as $(ij;kl;mn),$ 
where $\{i,j,k,l,m,n\}=\{1,\dots,6\}.$ Then 
it turns out that the diagram (\ref{diagram}) is $S_6$-equivariant. 

An explicit expression of $\psi^{-1}$ is given in \cite{Gon}.  
We want to know the combinatorial structure of $\psi^{-1}$ 
in order to study the inverse of the period map from a family of smooth cubic 
surfaces to the $4$-dimensional complex ball $\B^4$ in \cite{ACT}. 

For a $2$-dimensional subfamily of ours defined by $\l_5=\l_6,$ 
the period map and its inverse are studied in \cite{Pic} and \cite{Shi}.

\bigskip\noindent

{\bf Acknowledgments.} 
The author is grateful to Professors H. Shiga and 
M. Yoshida for their valuable discussions.

\section{Configuration space $\L$ of six points on $\Ps^1$}
Let $M(m,n)$ be the set of complex $(m\times n)$ matrices.  
We define the configuration space $\L$ of ordered six distinct points on 
the complex projective line $\Ps^1$ as 
$$\L=GL_2(\C)\backslash M'(2,6)/(\C^*)^6,$$
where 
$$M'(2,6)=\{\l=(\l_{ij})\in M(2,6)\mid \l\la kl\ra=
\left|\matrix{\l_{1k}&\l_{1l}\cr \l_{2k}&\l_{2l}\cr}\right|
\not= 0 \ (1\le k\not= l\le 6)\},$$
and $GL_2(\C)$ and $(\C^*)^6$ 
(regarding as the group of $(6\times 6)$ diagonal matrices)  
act naturally on $M'(2,6)$ from the left and right, respectively.
Note that we regard the column vectors of $\l\in M'(2,6)$ as the homogeneous 
coordinates of six points on $\Ps^1$ and the action of $GL_2(\C)$ as
the projective transformation.  Six distinct points $\l_1,\dots,\l_6$ on $\C$ 
are expressed by an element of $\L$ by $(2\times6)$ matrix 
$$\l=\pmatrix{  1 &  1 & 1  & 1  & 1  &  1 \cr 
              \l_1&\l_2&\l_3&\l_4&\l_5&\l_6\cr}.$$
By normalizing $(\l_1,\l_2,\l_3)$ as $(\infty,0,1),$ matrices of the form 
$$\pmatrix{ 0 & 1 & 1 & 1  & 1  &  1 \cr 
            1 & 0 & 1 &\ell_1&\ell_2&\ell_3\cr},\quad \ell_i\not= 0,1,\ell_j\ 
(1\le i<j\le 3)$$
represent $\L.$

We define a map  $\iota$ from $\L$ to the $14$-dimensional 
projective space $\Ps^{14}$ by 
$$\iota:\L\ni \l\mapsto 
[\dots,y_{\la ij; kl;mn\ra},\dots]=
[\dots,\l\la ij\ra \l\la kl\ra \l\la mn\ra,\dots]\in \Ps^{14},$$
where $\l$ is a $(2\times 6)$ matrix represent of an element of $\L$  
and projective coordinates of $\Ps^{14}$ are labeled by 
$I=\la ij;kl;mn\ra$ $(\{i,j,k,l,m,n\}=\{1,\dots,6\},\ i<j,\ k<l,\ m<n).$ 
Since the image $\iota(\l)$ is invariant under the actions $GL_2(\C)$ 
and $(\C^*)^6,$ this map is well defined. 
We use the following convention 
$$y_{\la ij; kl;mn\ra}=y_{\la kl; ij;mn\ra}=y_{\la ij; mn;kl\ra}
=-y_{\la ji; kl;mn\ra}.$$
The image $\iota(\L)$ is studied in \cite{Yo2}, 
it is described as the following.
\begin{fact}
The closure $Y=\overline{\iota(\L)}$ of $\iota(\L)$  is 
a subvariety of $\Ps^{14}$ defined by the linear and cubic equations
$$y_{\la ij; kl;mn\ra}-y_{\la ij; km;ln\ra}+y_{\la ij; kn;lm\ra}=0$$
$$y_{\la ij; kl;mn\ra}y_{\la ik; jn;lm\ra}y_{\la im; jl;kn\ra}=
y_{\la ij; kn;lm\ra}y_{\la ik; jl;mn\ra}y_{\la im; jn;kl\ra}.$$
\comment{
The variety $Y$ admits the following stratification. 
Let us define subvariety of $Y:$
\begin{eqnarray*}
Y'&=& \{y\in \Ps^{14}\mid   y_I\not= 0,\quad {\rm for\ any\ } I\},\ 
{\rm open\ in\ }Y,\\
Y^{ij}&=& \{y\in \Ps^{14}\mid  y_I= 0 {\rm\ for\ }  
(ij)\subset I,\ y_I\not= 0 {\rm\ otherwise}
\},\\
Y^{ij;kl}&=& 
\{y\in \Ps^{14}\mid  y_I=0 {\rm\ for\ } (ij) {\rm\ or\ } (kl) \subset I,\ 
y_I\not= 0 {\rm\ otherwise}
\},\\
Y^{ij;kl;mn}&=& 
\left\{y\in \Ps^{14}\mid
\matrix{  y_I=0 {\rm\ for\ } (ij) {\rm\ or\ } (kl) {\rm\ or\ } (mn) 
\subset I,\cr  y_I\not= 0 {\rm\ otherwise} }
\right\},\\
Y^{ijk;lmn}&=& 
\left\{y\in \Ps^{14}\mid
\matrix{  y_{\la pq;rs;tu\ra}\not= 0 {\rm \ for\ }
\{pq\},\{rs\},\{tu\}\not\subset \{i,j,k\}\cr {\rm \ or\ }\{l,m,n\},\quad 
y_I=0 {\rm\ otherwise} }\right\},
\end{eqnarray*}
where $(pq)\subset I=\la ij;kl;mn\ra$ means that 
one of $ij,$ $kl$ and $mn$ is $pq.$
The domain $Y'$ of $Y$ is a $3$-dimensional affine variety isomorphic to 
$\L;$  each $Y^{ij}$ is $2$-dimensional affine variety;
each $Y^{ij;kl}$ is $1$-dimensional affine variety;
each $Y^{ij;kl;mn}$ is a point; each $Y^{ijk;lmn}$ is a point.  
The variety $Y$ is irreducible and admits the following stratification:
$$Y=Y'\cup (\cap Y^{ij})\cup(\cup Y^{ij;kl})\cup(\cup Y^{ij;kl;mn})
\cup(\cup Y^{ijk;lmn}).$$
The variety $Y$ has singularities only at ten points $\cup Y^{ijk;lmn}.$
Let $\hat {\ } $ denote the closure in $Y;$ then
\begin{eqnarray*}
\hat Y^{ij}&=& Y^{ij}\cup(\cup \hat Y^{ij;kl})\simeq \Ps^2\\
\hat Y^{ij;kl}&=& Y^{ij;kl}\cup Y^{ij;kl;mn}
\cup Y^{ijm;kln}\cup Y^{ijn;klm}\simeq \Ps^1
\end{eqnarray*}
}
\end{fact}
We define $\hat\L$ as the compactification of $\L$ isomorphic to $Y.$

\section{Period matrix of $C$}
Let $C=C(\l)$ be the triple covering of $\Ps^1$ branching at six distinct 
points $\l_i'$s: 
$$C(\l): \ w^3=\prod_{i=1}^6 (z-\l_i);$$ 
this curve is of genus $4.$
Let $\rho$ be the automorphism of $C$ defined by 
$$\rho:C\ni (z,w) \mapsto (z,\w w)\in C,$$
where $\w={-1+\sqrt{-3}\over 2}.$
We give a basis of the vector space of holomorphic $1$-forms on 
$C$ as follows 
\begin{equation}
\label{cohom}
\f_1={dz\over w},\quad \f_2={dz\over w^2},\quad
\f_3={zdz\over w^2},\quad \f_4={z^2dz\over w^2}.
\end{equation}
For a fixed $\l$ such that $\l_i\in \R,\ \l_1<\dots<\l_6,$
we take a symplectic basis $\{A_1,\dots,A_4,B_1,\dots,B_4\}$ 
of $H_1(C,\Z)$ (i.e., $A_i\cdot A_j=B_i\cdot B_j=0,\ 
B_i\cdot A_j=\delta_{ij}$) such that 
\begin{equation}
\label{hom}
\rho(B_i)=A_i\ (i=1,2,3),\quad \rho(B_4)=-A_4,
\end{equation}
see Figure 1.

\begin{figure}[htb]
\begin{picture}(130,160)
\label{homolgy}
\thicklines
\multiput(-5,15)(0,50){3}{\line(1,0){130}}
\multiput(-5,55)(0,50){3}{\line(1,0){130}}
\multiput(-5,15)(0,50){3}{\line(0,1){40}}
\multiput(125,15)(0,50){3}{\line(0,1){40}}
\multiput(10,20)(0,50){3}{\line(1,0){100}}
\multiput(10,20)(20,0){6}{\line(0,1){20}}
\multiput(10,70)(20,0){6}{\line(0,1){20}}
\multiput(10,120)(20,0){6}{\line(0,1){20}}
\thinlines
\multiput(5,35)(40,0){3}{\line(0,1){10}}
\multiput(5,35)(40,0){3}{\line(1,0){5}}
\multiput(25,35)(40,0){3}{\line(1,0){5}}
\multiput(5,45)(40,0){3}{\line(1,0){10}}
\multiput(15,45)(40,0){3}{\vector(1,-1){10}}
\multiput(35,85)(40,0){3}{\line(0,1){10}}
\multiput(10,85)(40,0){3}{\line(1,0){5}}
\multiput(30,85)(40,0){3}{\line(1,0){5}}
\multiput(25,95)(40,0){3}{\line(1,0){10}}
\multiput(25,95)(40,0){3}{\vector(-1,-1){10}}
\put(30,30){\line(1,0){5}}
\put(35,30){\line(1,2){10}}
\put(45,50){\vector(1,0){40}}
\put(85,50){\line(1,0){35}}
\put(120,30){\line(0,1){20}}
\put(110,30){\line(1,0){10}}
\put(0,80){\line(1,0){10}}
\put(0,80){\line(0,1){20}}
\put(90,80){\line(-1,0){5}}
\put(85,80){\line(-1,2){10}}
\put(75,100){\vector(-1,0){40}}
\put(0,100){\line(1,0){35}}
\put(10,130){\vector(1,0){20}}
\put(110,130){\vector(-1,0){20}}
\textcolor{blue}{
\multiput(35,33)(40,0){3}{\line(0,1){10}}
\multiput(10,33)(40,0){3}{\line(1,0){5}}
\multiput(30,33)(40,0){3}{\line(1,0){5}}
\multiput(25,43)(40,0){3}{\line(1,0){10}}
\multiput(25,43)(40,0){3}{\vector(-1,-1){10}}
\multiput(5,133)(40,0){3}{\line(0,1){10}}
\multiput(5,133)(40,0){3}{\line(1,0){5}}
\multiput(25,133)(40,0){3}{\line(1,0){5}}
\multiput(5,143)(40,0){3}{\line(1,0){10}}
\multiput(15,143)(40,0){3}{\vector(1,-1){10}}
\put(0,28){\line(1,0){10}}
\put(0,28){\line(0,1){20}}
\put(90,28){\line(-1,0){5}}
\put(85,28){\line(-1,2){10}}
\put(40,48){\line(1,0){35}}
\put(0,48){\vector(1,0){40}}
\put(30,78){\vector(-1,0){20}}
\put(90,78){\vector(1,0){20}}
\put(30,128){\line(1,0){5}}
\put(35,128){\line(1,2){10}}
\put(45,148){\line(1,0){40}}
\put(120,148){\vector(-1,0){35}}
\put(120,128){\line(0,1){20}}
\put(110,128){\line(1,0){10}}
}
\put(20,0){Figure 1. \quad Our basis of $H_1(C,\Z)$}
\put(45,10){the first $z$-space}
\put(45,60){the second $z$-space}
\put(45,110){the third $z$-space}
\multiput(11,40)(0,50){3}{$\l_1$}
\multiput(31,40)(0,50){3}{$\l_2$}
\multiput(51,40)(0,50){3}{$\l_3$}
\multiput(71,40)(0,50){3}{$\l_4$}
\multiput(91,40)(0,50){3}{$\l_5$}
\multiput(111,40)(0,50){3}{$\l_6$}
\put(17,43){$A_1$}
\put(57,43){$A_2$}
\put(97,43){$A_3$}
\put(80,47){$A_4$}
\textcolor{blue}{
\put(17,33){$B_1$}
\put(57,33){$B_2$}
\put(97,33){$B_3$}
\put(32,49){$B_4$}
}
\put(17,85){$A_1$}
\put(57,85){$A_2$}
\put(97,85){$A_3$}
\put(40,97){$A_4$}
\textcolor{blue}{
\multiput(20,74)(80,0){2}{$B_4$}
}
\multiput(20,127)(80,0){2}{$A_4$}
\textcolor{blue}{
\put(17,142){$B_1$}
\put(57,142){$B_2$}
\put(97,142){$B_3$}
\put(80,149){$B_4$}
}
\end{picture}
\end{figure}
\comment{
\begin{figure}[htb]
\begin{picture}(130,160)
\label{homolgy}
\thicklines
\multiput(-5,15)(0,50){3}{\line(1,0){130}}
\multiput(-5,55)(0,50){3}{\line(1,0){130}}
\multiput(-5,15)(0,50){3}{\line(0,1){40}}
\multiput(125,15)(0,50){3}{\line(0,1){40}}
\multiput(10,20)(0,50){3}{\line(1,0){100}}
\multiput(10,20)(20,0){6}{\line(0,1){20}}
\multiput(10,70)(20,0){6}{\line(0,1){20}}
\multiput(10,120)(20,0){6}{\line(0,1){20}}
\thinlines
\multiput(5,35)(40,0){3}{\line(0,1){10}}
\multiput(5,35)(40,0){3}{\line(1,0){5}}
\multiput(25,35)(40,0){3}{\line(1,0){5}}
\multiput(5,45)(40,0){3}{\line(1,0){10}}
\multiput(15,45)(40,0){3}{\vector(1,-1){10}}
\multiput(35,85)(40,0){3}{\line(0,1){10}}
\multiput(10,85)(40,0){3}{\line(1,0){5}}
\multiput(30,85)(40,0){3}{\line(1,0){5}}
\multiput(25,95)(40,0){3}{\line(1,0){10}}
\multiput(25,95)(40,0){3}{\vector(-1,-1){10}}
\put(30,30){\line(1,0){5}}
\put(35,30){\line(1,2){10}}
\put(45,50){\vector(1,0){20}}
\put(65,50){\line(1,0){15}}
\put(80,30){\line(0,1){20}}
\put(70,30){\line(1,0){10}}
\put(0,80){\line(1,0){10}}
\put(0,80){\line(0,1){20}}
\put(50,80){\line(-1,0){5}}
\put(45,80){\line(-1,2){10}}
\put(35,100){\vector(-1,0){20}}
\put(0,100){\line(1,0){15}}
\put(10,130){\vector(1,0){20}}
\put(70,130){\vector(-1,0){20}}
\textcolor{blue}{
\multiput(35,33)(40,0){3}{\line(0,1){10}}
\multiput(10,33)(40,0){3}{\line(1,0){5}}
\multiput(30,33)(40,0){3}{\line(1,0){5}}
\multiput(25,43)(40,0){3}{\line(1,0){10}}
\multiput(25,43)(40,0){3}{\vector(-1,-1){10}}
\multiput(5,133)(40,0){3}{\line(0,1){10}}
\multiput(5,133)(40,0){3}{\line(1,0){5}}
\multiput(25,133)(40,0){3}{\line(1,0){5}}
\multiput(5,143)(40,0){3}{\line(1,0){10}}
\multiput(15,143)(40,0){3}{\vector(1,-1){10}}
\put(0,28){\line(1,0){10}}
\put(0,28){\line(0,1){20}}
\put(50,28){\line(-1,0){5}}
\put(45,28){\line(-1,2){10}}
\put(20,48){\line(1,0){15}}
\put(0,48){\vector(1,0){20}}
\put(30,78){\vector(-1,0){20}}
\put(50,78){\vector(1,0){20}}
\put(30,128){\line(1,0){5}}
\put(35,128){\line(1,2){10}}
\put(45,148){\line(1,0){20}}
\put(80,148){\vector(-1,0){15}}
\put(80,128){\line(0,1){20}}
\put(70,128){\line(1,0){10}}
}
\put(20,0){Figure 1. \quad Our basis of $H_1(C,\Z)$}
\put(45,10){the first $z$-space}
\put(45,60){the second $z$-space}
\put(45,110){the third $z$-space}
\multiput(11,40)(0,50){3}{$\l_1$}
\multiput(31,40)(0,50){3}{$\l_2$}
\multiput(51,40)(0,50){3}{$\l_3$}
\multiput(71,40)(0,50){3}{$\l_4$}
\multiput(91,40)(0,50){3}{$\l_5$}
\multiput(111,40)(0,50){3}{$\l_6$}
\put(17,43){$A_1$}
\put(57,43){$A_2$}
\put(97,43){$A_3$}
\put(80,47){$A_4$}
\textcolor{blue}{
\put(17,33){$B_1$}
\put(57,33){$B_2$}
\put(97,33){$B_3$}
\put(32,49){$B_4$}
}
\put(17,85){$A_1$}
\put(57,85){$A_2$}
\put(97,85){$A_3$}
\put(40,97){$A_4$}
\textcolor{blue}{
\multiput(20,74)(40,0){2}{$B_4$}
}
\multiput(20,127)(40,0){2}{$A_4$}
\textcolor{blue}{
\put(17,142){$B_1$}
\put(57,142){$B_2$}
\put(97,142){$B_3$}
\put(80,149){$B_4$}
}
\end{picture}
\end{figure}}

Put 
$${\int_{A_i}\f_j\choose\int_{B_i}\f_j }_{i,j}=
{\W_A\choose\W_B}.$$
Let $\f$ be the normalized basis of vector space of 
holomorphic $1$-forms so that $\W_B$ becomes $I_4.$ 
Note that the normalized period $\W=\W_A\W_B^{-1}$ 
belongs to the Siegel upper half space $\S^4$ of degree $4.$ 
The next proposition shows that $\W$ can be expressed in terms of 
$$x=\tr(x_1,\dots,x_4)=\tr\left(\int_{A_1}\f_1,\dots,\int_{A_4}\f_1\right).$$
\begin{proposition}
\label{permat}
We have 
\begin{eqnarray*}
\W&=&\w[I_4-(1-\w)(x\tr xH)/(\tr x Hx)]H=\w[H-(1-\w)(x\tr x)/(\tr x Hx)]\\
&=&\pmatrix{
\w&  &  &  \cr
  &\w&  &  \cr
  &  &\w&  \cr
  &  &  &-\w\cr}
\!-\! {\sqrt{-3}\over x_1^2\!+\! x_2^2\!+\! x_3^2\!-\! x_4^2}
\pmatrix{
x_1x_1& x_1x_2 &x_1x_3 &x_1x_4 \cr
x_2x_1& x_2x_2 &x_2x_3 &x_2x_4 \cr
x_3x_1& x_3x_2 &x_3x_3 &x_3x_4 \cr
x_4x_1& x_4x_2 &x_4x_3 &x_4x_4 \cr},
\end{eqnarray*}
where $H={\rm diag}(1,1,1,-1)$ and $\tr{\bar x} H x<0.$
\end{proposition}
{\itshape Proof.} Put $\W_A=({ x},{ b},{ c},{ d});$
by (\ref{cohom}) and (\ref{hom}), $\W_B$ can be expressed as 
$$\W_B=(\w H{ x},\w^2H{ b},\w^2H{ c},\w^2H{ d}) 
=\w^2 H\W_A+(\w-\w^2)H(x,O).$$
We have 
$$\W^{-1}=\W_B\W_A^{-1}=\w^2 H +(\w-\w^2)H(x,O)\W_A^{-1}.$$
Put 
$$\W_A^{-1}={\xi\choose \ast},\quad \xi=(\xi_1,\dots,\xi_4);$$
note that 
$$\xi x=\sum_{i=1}^4 \xi_ix_i=1.$$
We have 
\begin{equation}
\label{inv}
H(x,O)\W_A^{-1}=Hx\xi={1\over \xi x}\pmatrix{
x_1\xi_1& x_1\xi_2& x_1\xi_3 & x_1\xi_4\cr
x_2\xi_1& x_2\xi_2& x_2\xi_3 & x_2\xi_4\cr
x_3\xi_1& x_3\xi_2& x_3\xi_3 & x_3\xi_4\cr
-x_4\xi_1& -x_4\xi_2& -x_4\xi_3& -x_4\xi_4\cr},
\end{equation}
which must be symmetric. Thus we have 
$$x_i\xi_j=x_j\xi_i\ (1\le i<j\le 3),\quad x_i\xi_4=-x_4\xi_i\ (i=1,2,3).$$
By eliminating $\xi_i$ in (\ref{inv}),   we have
$$H(x,O)\W_A^{-1}=(Hx\tr xH)/(\tr x Hx).$$
Then 
$$\W^{-1}=\w^2H[I_4-(1-\w^2)(x\tr xH)/(\tr x Hx)].$$
It is easy to see that 
$$[I_4-(1-\w^2)(x\tr xH)/(\tr x Hx)]^{-1}=I_4-(1-\w)(x\tr xH)/(\tr x Hx),$$
we have
$$
\W=\w[I_4-(1-\w)(x\tr xH)/(\tr x Hx)]H.$$

The imaginary part of $\W$ is ${\sqrt3\over 2}$ times 
\begin{equation}
\label{impart}
H-x\tr x/(\tr x H x)-\bar x\tr \bar x/(\tr\bar x H\bar x),
\end{equation}
which must be positive definite.
If $x_4=0$ then the $(4,4)$ component of (\ref{impart}) is $-1,$
which implies that (\ref{impart}) can not be positive definite. Thus  
we have $x_4\not= 0.$ 
Put 
$$\eta=\pmatrix{x_4& 0 & 0 \cr
                 0 &x_4& 0 \cr
                 0 & 0 &x_4\cr
                x_1&x_2&x_3\cr};$$
note that $(\eta,x)\in GL_4(\C)$ and that $\tr x H\eta=(0,0,0).$
We have 
$$\tr\overline{(\eta,x) H}\left(H-{x\tr x\over \tr x H x}-
{\bar x\tr \bar x\over \tr\bar x H\bar x}\right)H(\eta,x)=
\pmatrix{
\tr\bar \eta H\eta & 0 \cr
 \tr 0 & -\tr \bar x H x\cr}.$$
\comment{
\begin{eqnarray*}
& &\tr\overline{(\eta,x) H}\left(H-{x\tr x\over \tr x H x}-
{\bar x\tr \bar x\over \tr\bar x H\bar x}\right)H(\eta,x)\\
&=&
\pmatrix{\tr\bar \eta\cr \tr \bar x} H (\eta,x) 
-\pmatrix{\tr\bar \eta Hx \cr \tr \bar xHx} {\tr xH(\eta,x)\over \tr x H x}
-{\pmatrix{\tr\bar \eta\cr \tr \bar x}H\bar x\over \tr\bar x H\bar x}
(\tr \bar xH\eta,\tr \bar xHx)\\
&=&\pmatrix{\tr\bar\eta H\eta &\tr\bar\eta H x \cr
          \tr\bar x H\eta & \tr\bar xHx}-\pmatrix{\tr\bar\eta H x\cr \bar x Hx}
(0,0,0,1)-\tr(0,0,0,1)(\tr\bar xH\eta,\tr\bar xHx)\\
&=&\pmatrix{
\tr\bar \eta H\eta & 0 \cr
 \tr 0 & -\tr \bar x H x\cr}.
\end{eqnarray*}}
If 
$$ -\bar x H x=-|x_1|^2-|x_2|^2-|x_3|^2+|x_4|^2>0$$
then the $3\times3$ matrix
$$\tr\bar \eta H\eta=
|x_4|^2I_3-\pmatrix{\bar x_1\cr\bar x_2\cr\bar x_3}(x_1,x_2,x_3)$$
is positive definite. 
Hence the matrix (\ref{impart}) is positive definite if and only if 
$$ ^t\bar x H x=|x_1|^2+|x_2|^2+|x_3|^2-|x_4|^2<0.\qedd$$

We embedded the domain 
$\B^3=\{x\in \Ps^3\mid \tr\bar x H x<0\}$ 
in $\S^4$ by the map 
$$\jmath:\B^3\ni x\mapsto \W=\w[I_4-(1-\w)(x\tr xH)/(\tr x Hx)]H\in \S^4.$$


\section{Monodromy }
Let $(\l_1,\dots,\l_6)$ vary as an element in $\L,$ 
we have two multi-valued map 
$$\matrix{\psi:& \L& \to& \B^3\cr
               & \l& \mapsto &  
x=\tr\left(\int_{A_1}\f_1,\dots,\int_{A_4}\f_1\right),\cr
 & & \cr
\tilde\psi=\jmath\circ\psi:&  \L & \to &\S^4\cr
                     & \l  &\mapsto&\W=\jmath(\psi(\l)).\cr}$$
We call them period maps. 
The map $\psi$ and its monodromy group were studied in \cite{DM}, \cite{Ter}, 
\cite{Yo1} and \cite{Yo2}, the results are as follows.
\begin{fact}
\label{ball}
The image of $\psi$ is open dense in $\B^3.$ 
The monodromy group of $\psi$ is conjugate to the congruence subgroup
$$\G(1-\w)=\{g\in \G\mid g\equiv I_4\ \mod \ (1-\w)\}$$
of the modular group 
$$\G=\{g\in GL_4(\Z[\w])\mid \tr\bar g H g=H\}.$$
The Satake compactification $\hat \B^3/\G(1-\w)$ of  
$\B^3/\G(1-\w)$ is isomorphic to $Y.$
\end{fact}  

For a column vector $v\in \C^4$ such that $\tr\bar v H v\not= 0,$ 
we define reflections $R^\w(v)$ and $R^\z(v)$  with  root $v$ and 
exponent $\w$ and $\z=-\w^2,$ respectively, as
\begin{eqnarray*}
R^\w(v)&=&I_4-(1-\w)v(\tr\bar v H v)^{-1}\tr\bar v H,\\
R^\z(v)&=&I_4-(1-\z)v(\tr\bar v H v)^{-1}\tr\bar v H.
\end{eqnarray*}
It is shown in \cite{All} that $\G(1-\w)$  can be generated by fifteen 
reflections
$R^\w_{ij}=R^\w(v_{ij})$ $(1\le i<j\le 6)$ 
and that $\G$ by $-I_4$ and five reflections 
$R^\z_{i,i+1} =R^\z(v_{i,i+1})$ $(1\le i\le 5),$
where 
$$
\begin{array}{lll}
v_{12}= \tr(1,0,0,0), &v_{13}= \tr(-1,1,0,1), & 
v_{14}= \tr(-1,-\w^2,0,1),\\
v_{15}= \tr(\w^2,0,-\w^2,1), &v_{16}= \tr(\w^2,0,\w,1), 
& v_{23}= \tr(\w^2,1,0,1),\\
v_{24}= \tr(\w^2,-\w^2,0,1), &v_{25}= \tr(-\w,0,-\w^2,1), 
& v_{26}= \tr(-\w,0,\w,1),\\
v_{34}= \tr(0,1,0,0), &v_{35}= \tr(0,-\w,\w,1), 
& v_{36}= \tr(0,-\w,-1,1),\\
v_{45}= \tr(0,1,\w,1), &v_{46}= \tr(0,1,-1,1), & v_{56}= \tr(0,0,1,0).\\
\end{array}
$$
The reflections correspond to the following movements of $\l_i$'s. 
When $\l_i$ goes near to $\l_j$ in the upper half space and turns around 
$\l_j$ and returns, $x$ becomes $R^\w_{ij}x.$ 
When $\l_i$ and $\l_j$ are exchanged in the upper half space,  
$x$ becomes $R^\z_{ij}x.$ Since $R^\z_{i,i+1}$'s are 
representations of braids, they satisfy 
$$R^\z_{i-1,i}R^\z_{i,i+1}R^\z_{i-1,i}
=R^\z_{i,i+1}R^\z_{i-1,i}R^\z_{i,i+1}\quad (2\le i\le 5).$$

The embedding $\jmath$ induces the following homomorphism 
from $U(3,1;\C)$ to 
$$Sp(8,\R)=\left\{g\in GL_8(\R)\mid \tr g J g=J=
\pmatrix{ O & -I_4\cr I_4 & O}\right\}:$$ 
$$\tilde\jmath:U(3,1;\C)\ni 
P+\w Q \mapsto \pmatrix{
P & QH\cr
-HQ & H(P-Q)H\cr}\in Sp(8,\R),$$
where $P$ and $Q$ are real $4\times 4$ matrices. 
Note that 
\begin{eqnarray*}\tilde\jmath^{-1}:Sp(8,\R)&\supset&\tilde\jmath(U(3,1;\C))
\ni \pmatrix{A & B\cr
C& D\cr} \cr
&\mapsto& A+\w BH=(-HC+HDH)-\w HC \in U(3,1;\C).
\end{eqnarray*}

Let us express the images of $R^\w(v)$ and $R^\z(v)$ 
under the map $\tilde\jmath.$ 
The image of $\w I_4$ under $\tilde\jmath$ is given by 
$$W=\pmatrix{O & H\cr -H & -I_4\cr}\in Sp(8,\Z).$$ 

For a column vector $v=a+\w b$ $(a,b\in \R^4),$
define column vectors $v_1={a \choose  {-Hb}}$ and $v_2=Wv_1$ and 
form a $(8\times 2)$ matrix $V=(v_1,v_2).$ 
Straightforward calculation shows the following. 
\begin{proposition}
\label{Grass}
If $\tr \bar v H v\not=0,$ then $\tilde\jmath(R^{\w}(v))=\tilde R^{\w}(v)$ and 
$\tilde\jmath(R^\z(v))=\tilde R^\z(v)$ are given by
$$I_{8}-(I_{8}-W)V(\tr VJV)^{-1}\tr VJ,\quad 
I_{8}-(I_{8}+W^2)V(\tr VJV)^{-1}\tr VJ,$$
respectively. 
\end{proposition}

Systems of generators of $\tilde\G(1-\w)=\tilde\jmath(\G(1-\w))$ and 
$\tilde\G=\tilde\jmath(\G)$ are given by 
$\tilde R^\w_{ij}$'s and $\tilde R^\z_{i,i+1}$'s.


\section{Riemann theta constants}
The Riemann theta function 
$$\vartheta(z,\tau)=
\sum_{n=(n_1,\dots,n_r)\in \Z^r}\exp[\pi\sqrt{-1}(n \tau\tr  n+2n \tr z)]$$
is holomorphic on $\C^r\times \S^r$ and satisfies
$$\vartheta(z+p,\tau)=\vartheta(z,\tau),\quad 
\vartheta(z+p\tau,\tau)=\exp[-\pi\sqrt{-1}(p\tau \tr p+2z\tr p )]
\vartheta(z,\tau),$$
where $\S^r$ 
is the Siegel upper half space of degree $r$ and $p\in \Z^r.$
It is well known that for $(z,\tau)\in \C\times \H,$
$\vartheta(z,\tau)=0$ if and only if $z={1+\tau\over 2}+p+q\tau\ (p,q\in\Z).$ 

The theta function $\vartheta_{a,b}(z,\tau)$ 
with characteristics $a,b$ is defined by 
\begin{eqnarray}
\label{thetaconst}
\vartheta_{a,b}(z,\tau)&=&\exp[\pi\sqrt{-1}(a\tau \tr a+2a\tr (z+b))]
\vartheta(z+a\tau+b,\tau)\\
&=&
\sum_{n\in \Z^n}\exp[\pi\sqrt{-1}((n+a) \tau \tr (n+a)+2(n+a)\tr(z+b))],
\nonumber
\end{eqnarray}
where $a,b\in \Q^r.$ Note that 
\begin{equation}
\label{Z-shift}
\vartheta_{-a,-b}(z,\tau)=\vartheta_{a,b}(-z,\tau),\quad 
\vartheta_{a+p,b+q}(z,\tau)=\exp(2\pi\sqrt{-1}a\tr q)\vartheta_{a,b}(z,\tau).
\end{equation}
The function $\vartheta_{a,b}(\tau)= \vartheta_{a,b}(0,\tau)$ of $\tau$ 
is called the theta constant with characteristics $a,b.$  
If $\tau$ is diagonal, then this function becomes the product of 
Jacobi's theta constants: 
$$\vartheta_{a,b}(\tau)=\prod_{i=1}^r\vartheta_{a_i,b_i}(\tau_i),$$
where 
$$a=(a_1,\dots,a_r),\ b=(b_1,\dots,b_r),\ \tau={\rm diag}(\tau_1,\dots,\tau_r).
$$

The following transformation formula can be found in \cite{Igu} p.176.
\begin{fact}
\label{trans}
For any $g= \pmatrix{A &B \cr C& D\cr}\in Sp(2r,\Z)$ and $(a,b)\in \Q^{2r},$ 
we put 
\begin{eqnarray*}
g\cdot(a,b)&=&(a,b)g^{-1}+{1\over 2}({\rm dv}(C\tr D),{\rm dv}(A\tr B))\\
\phi_{(a,b)}(g)&=&-{1\over 2}(a\tr D B\tr a-2a\tr B C\tr b+b\tr C A\tr b)\\
& & +{1\over 2}(a\tr D-b\tr C)\tr({\rm dv}(A\tr B)),
\end{eqnarray*}
where ${\rm dv}(A)$ is the row vector consisting of 
the diagonal components of $A.$
Then for every $g\in Sp(2r,\Z),$ we have 
\begin{eqnarray*}
& &\vartheta_{g\cdot (a,b)}((A\tau +B)(C\tau +D)^{-1})\\
&=&\kappa(g)\exp(2\pi\sqrt{-1}\phi_{(a,b)}(g))\det(C\tau+D)^{1\over 2}
\vartheta_{(a,b)}(\tau),
\end{eqnarray*}
in which $\kappa(g)^2$ is a $4$-th root of $1$ depending only on $g.$ 
\end{fact}

\begin{proposition}
\label{invchar}
There are $81=3^4$ theta characteristics 
$$(a,b)=(a_1,\dots,a_4,b_1,\dots,b_4)$$ 
such that 
$$g\cdot(a,b)\equiv (a,b)\ \mod \Z^8$$
for any $g\in \tilde\G(1-\w)\subset Sp(8,\Z);$
they are given by 
\begin{equation}
\label{cond}
b=-aH, \ a_i\in \{{1\over 6},{3\over 6},{5\over 6}\}\ (i=1,\dots,4).
\end{equation}
\end{proposition}
{\itshape Proof.} 
Since 
$$W\cdot(a,b)=(-a+bH,-aH)+{1\over 2}(1,1,1,-1,0,0,0,0),$$
we have 
$$-aH \equiv b,\ -2a+{1\over 2}(1,1,1,-1)\equiv a\ \mod \Z^4.$$
Thus we have the condition (\ref{cond}). 
It is easy to check such theta characteristics are invariant under the action 
on $15$ reflections $\tilde R^\w_{ij}.$ \qed
We label the $81$ characteristics $a$'s by combinatorics of six letters;  
they are classified to $4$ classes. 
The list of the correspondence between the label of $a$ and $6a$ is as follows:
{\footnotesize
$$\matrix{
(12;\!34;\!56)\!\leftrightarrow\!\pm\!({3},{3},{3},{-1})& 
(12;\!35;\!46)\!\leftrightarrow\!\pm\!({3},{1},{1},{-3})&
(12;\!36;\!45)\!\leftrightarrow\!\pm\!({3},{1},{-1},{-3})\cr
& & \cr
(13;\!24;\!56)\!\leftrightarrow\!\pm\!({1},{1},{3},{-3})& 
(13;\!25;\!46)\!\leftrightarrow\!\pm\!({1},{-1},{1},{-1})&
(13;\!26;\!45)\!\leftrightarrow\!\pm\!({-1},{1},{1},{1})\cr
& & \cr
(14;\!23;\!56)\!\leftrightarrow\!\pm\!({1},{-1},{3},{-3})& 
(14;\!25;\!36)\!\leftrightarrow\!\pm\!({1},{1},{1},{-1})&
(14;\!26;\!35)\!\leftrightarrow\!\pm\!({1},{1},{-1},{-1})\cr
& & \cr
(15;\!23;\!46)\!\leftrightarrow\!\pm\!({1},{1},{-1},{1})&
(15;\!24;\!36)\!\leftrightarrow\!\pm\!({1},{-1},{-1},{1})&
(15;\!26;\!34)\!\leftrightarrow\!\pm\!({1},{3},{1},{-3})\cr
& & \cr
(16;\!23;\!45)\!\leftrightarrow\!\pm\!({1},{1},{1},{1})& 
(16;\!24;\!35)\!\leftrightarrow\!\pm\!({1},{-1},{1},{1})&
(16;\!25;\!34)\!\leftrightarrow\!\pm\!({1},{3},{-1},{-3})\cr
}$$}

$$\matrix{
(1^22)\!\leftrightarrow\!({1},{3},{3},{3})& 
(1^23)\!\leftrightarrow\!({5},{1},{3},{5})&
(1^24)\!\leftrightarrow\!({5},{5},{3},{5})\cr
& & \cr
(1^25)\!\leftrightarrow\!({5},{3},{1},{1})& 
(1^26)\!\leftrightarrow\!({5},{3},{5},{1})&
(2^23)\!\leftrightarrow\!({1},{1},{3},{5})\cr
& & \cr
(2^24)\!\leftrightarrow\!({1},{5},{3},{5})& 
(2^25)\!\leftrightarrow\!({1},{3},{1},{1})&
(2^26)\!\leftrightarrow\!({1},{3},{5},{1})\cr
& & \cr
(3^24)\!\leftrightarrow\!({3},{1},{3},{3})& 
(3^25)\!\leftrightarrow\!({3},{5},{1},{5})&
(3^26)\!\leftrightarrow\!({3},{5},{5},{5})\cr
& & \cr
(4^25)\!\leftrightarrow\!({3},{1},{1},{5})& 
(4^26)\!\leftrightarrow\!({3},{1},{5},{5})&
(5^26)\!\leftrightarrow\!({3},{3},{1},{3})\cr
}$$
$$(ij^2)\leftrightarrow -a {\rm \ for\ }(i^2j)\quad 1\le i<j\le 6,$$

$$\matrix{
(123)\!\leftrightarrow\!({3},{1},{3},{5})& 
(124)\!\leftrightarrow\!({3},{5},{3},{5})&
(125)\!\leftrightarrow\!({3},{3},{1},{1})\cr
& & \cr
(126)\!\leftrightarrow\!({3},{3},{5},{1})& 
(134)\!\leftrightarrow\!({1},{3},{3},{1})&
(135)\!\leftrightarrow\!({1},{1},{1},{3})\cr
& & \cr
(136)\!\leftrightarrow\!({1},{1},{5},{3})& 
(145)\!\leftrightarrow\!({1},{5},{1},{3})&
(146)\!\leftrightarrow\!({1},{5},{5},{3})\cr
& & \cr
(156)\!\leftrightarrow\!({1},{3},{3},{5})& 
(lmn)\leftrightarrow-a{\rm \ for\ }(ijk) & \{i,j,k,l,m,n\}=\{1,\dots,6\}
\cr
}$$

$$(123456)\!\leftrightarrow\!({3},{3},{3},{3}).$$
The first class is characterized by $(6a) H \tr(6a)\equiv 2 \mod 24$ and 
the characteristics $(a,-aH)$ with label $(ij;kl;mn)$ is 
invariant under the actions $\tilde R^\z_{ij},$ $\tilde R^\z_{kl}$ and 
$\tilde R^\z_{mn};$ 
the second class is characterized by $(6a) H \tr(6a)\equiv 10 \mod 24$ and 
the characteristics $(a,-aH)$ with label $(i^2j)$ is 
invariant under the actions $\tilde R^\z_{kl}$ 
$(\{i,j\}\cap\{k,l\}=\emptyset)$ 
and $\tilde R^\z_{ij}\cdot(a,-aH)$ is $(-a,aH)$ with label $(ij^2);$ 
the third class is characterized by $(6a) H \tr(6a)\equiv 18 \mod 24$ and 
the characteristics $(a,-aH)$ with label $(ijk)$ is 
invariant under the actions $\tilde R^\z_{lm}$ $(\{i,j,k\}\cap\{l,m\}=
\emptyset\ {\rm or\ } \{l,m\}).$ 

We denote $\vartheta_{a,-aH}(\W)$ by $\vartheta_{[6a]}(\W)$ or 
$\vartheta(ij;kl;mn),$ 
$\vartheta(i^2j),$ $\vartheta(ijk)$ and $\vartheta(123456)$ for 
corresponding characteristics $a.$
Note that for $p,q\in \Z^4,$
\begin{eqnarray*}
& & \vartheta(a(\W\!-\! H)+p\W+q,\W)\\
&=&\exp[-\pi\sqrt{-1}(p\W\tr p\!+\!2p(\W\!-\! H)\tr a)]
\vartheta(a(\W\!-\! H),\W)\\
&=&\exp[-\pi\sqrt{-\!1}(p\W\tr p\!+\!2p(\W\!-\! H)\tr a
\!+\! a\W\tr a\!-\! 2aH\tr a)]\vartheta_{a,-aH}(\W)\\
&=&\exp[2\pi\sqrt{-1}(a+p)H\tr(a+p)]
\exp[-\pi\sqrt{-1}(a+p)\W\tr(a+p)]\vartheta_{[a]}(\W).
\end{eqnarray*}
\begin{proposition}
\label{survive}
The theta constants $\vartheta(i^2j),$ $\vartheta(ijk)$ and 
$\vartheta(123456)$ are identically zero on $\jmath(\B^3).$ 
The theta constants $\vartheta(ij;kl;mn)$ are not 
identically zero on $\jmath(\B^3).$
\end{proposition}
{\itshape Proof.}
We apply Fact \ref{trans} for 
$$\tau=\W=\jmath(x),\quad 
g=W=\pmatrix{0 & H\cr -H& -I_4\cr},\quad (a,b)=(a,-aH).$$
Note that 
$$W\cdot \W=\W,\quad W\cdot(a,-aH)=
\left(a-(3a-{1\over 2}{\rm diag}(H)),-aH\right)$$
and that 
$$\phi_{(a,-aH)}(W)={3\over 2}a H\tr  a={1\over 24}(6a) H\tr(6a) ,\quad 
\det(C\W+D)=\w.$$
Since $\kappa(W)$ is an $8$-th root of $1,$ 
the sufficient condition for 
\begin{equation}
\label{8th}
\kappa(W)\exp(2\pi\sqrt{-1}\phi_{(a,b)}(W))\det(C\W+D)^{1\over 2}=1
\end{equation}
is $(6a)H\tr(6a)\equiv 2$ $ {\rm mod}\ 24.$
If $(6a)H\tr(6a)\not \equiv 2$ $ {\rm mod}\ 24,$ then $\vartheta_{a,-aH}(\W)$ 
vanishes. Thus the theta constants $\vartheta(i^2j),$ $\vartheta(ijk)$ and 
$\vartheta(123456)$ are identically zero on $\jmath(\B^3).$

For $a=({1\over 6},\dots,{1\over 6})$ and $x=(0,0,0,1),$ 
$\vartheta_{a,-aH}(\W)$ reduces to
$$\vartheta_{({1\over 6},{-1\over 6})}(\w)^3
\vartheta_{({1\over 6},{1\over 6})}(-\w^2),$$
which does not vanish.  Hence $\vartheta(ij;kl;mn)$'s survive. 
Note that $\kappa(W)^2=-1$  by (\ref{8th}).
\qed 

\begin{proposition}
\label{G-act}
We have 
\begin{eqnarray*}
\vartheta(i,i+1;kl;mn)(\tilde R^\z_{i,i+1}\cdot \jmath(x))^3
&=&-\chi(\tilde R^\z_{i,i+1})\vartheta(i,i+1;kl;mn)(\jmath(x))^3,\\
\vartheta(ik;i+1,l;mn)(\tilde R^\z_{i,i+1}\cdot \jmath(x))^3
&=&\chi(\tilde R^\z_{i,i+1})\vartheta(il;i+1,k;mn)(\jmath(x))^3,
\end{eqnarray*}
where
$$\chi(\tilde R^\z_{i,i+1})=
\left({\tr ( R^\z_{i,i+1} x) H ( R^\z_{i,i+1}x)\over  \tr xHx}\right)^{3/2},$$
which takes $1$ on the mirror of $R^\z_{i,i+1}.$
\end{proposition}
{\itshape Proof.} For $\tilde R^\z_{i,i+1}=\pmatrix{ A & B\cr C& D\cr},$
straightforward calculation shows 
$$\det(C\jmath(x)+D)={\tr ( R^\z_{i,i+1} x) H ( R^\z_{i,i+1}x)\over
 \det( R^\z_{i,i+1}) \tr xHx}
={\tr ( R^\z_{i,i+1} x) H ( R^\z_{i,i+1}x)\over -\w^2 \tr xHx}.$$
By computing $\phi_{a,b}(\tilde R^\z_{i,i+1})$ in Fact \ref{trans} and 
using (\ref{Z-shift}), 
we have 
\begin{eqnarray*}
\vartheta(i,i+1;kl;mn)(\tilde R^\z_{i,i+1}\cdot \jmath(x))^3
&=&-c\chi(\tilde R^\z_{i,i+1})\vartheta(i,i+1;kl;mn)(\jmath(x))^3,\\
\vartheta(ik;i+1,l;mn)(\tilde R^\w_{i,i+1}\cdot \jmath(x))^3
&=&c\chi(\tilde R^\z_{i,i+1})\vartheta(il;i+1,k;mn)(\jmath(x))^3,
\end{eqnarray*}
where $c$ is a certain constant depending only on $\tilde R^\z_{i,i+1}.$ 
If we restrict $\jmath(x)$ on the mirror of $R^\z_{i,i+1},$ 
we have 
$$\tilde R^\z_{i,i+1}\cdot \jmath(x)=\jmath(x),\quad \chi(\tilde R^\z_{i,i+1})=
\left({\tr (R^\z_{i,i+1} x) H ( R^\z_{i,i+1}x)\over \tr xHx}\right)^{3/2}=1.$$
Since $\vartheta(i,k;i+1,l;mn)=\vartheta(i,l;i+1,k;mn)$ 
on the mirror of $R^\z_{i,i+1}$ and it does not vanish,  
the constant $c$ must be $1.$ \qed

Since $\tilde R^\z_{pq}$ can be expressed in terms of 
$\tilde R^\z_{i,i+1}$ and $\tilde R^\w_{pq}=(\tilde R^\z_{pq})^2,$ 
we have the following two propositions.
\begin{proposition}
\label{G(1-w)-act}
We have 
$$
\vartheta(ij;kl;mn)(\tilde R^\w_{pq}\cdot \jmath(x))^3
=\chi(\tilde R^\w_{pq})\vartheta(ij;kl;mn)(\jmath(x))^3,
$$
where
$$\chi(\tilde R^\w_{pq})=
\left({\tr ( R^\w_{pq} x) H ( R^\w_{pq}x)\over  \tr xHx}\right)^{3/2},$$
which takes $1$ on the mirror of $R^\w_{pq}.$
\end{proposition}

\begin{proposition}
\label{vanish}
The function $\vartheta(ij;kl;mn)(\jmath(x))$ 
vanishes on the $\G(1-\w)$ orbits of 
the mirrors of $R^\w_{ij},$ $R^\w_{kl}$ and $R^\w_{mn}.$
\end{proposition}
{\itshape Proof.} By Proposition \ref{G-act}, when we restrict $\jmath(x)$ 
on the mirrors of $R^\w_{12},$ $R^\w_{34}$ and $R^\w_{56},$ we have
$$\vartheta(12;34;56)(\jmath(x))^3=-\vartheta(12;34;56)(\jmath(x))^3=0.$$ 
For the $\G(1-\w)$ orbits, use the previous proposition. 
In oder to show for general $\vartheta(ij;kl;mn)(\jmath(x))$'s,
use Proposition \ref{G-act}. \qed 


\section{The inverse of the period map} 
\begin{proposition}
\label{invper} Let $\W$ be the period matrix of 
$$C(\l):w^3=z(z-1)(z-\ell_1)(z-\ell_2)(z-\ell_3)$$
given in Proposition \ref{permat}. We have 
\begin{eqnarray}
\label{lambda4}
\ell_1&=&{\vartheta^3(13;24;56)(\W)\over \vartheta^3(14;23;56)(\W)},\\
\label{lambda5}
\ell_2&=&{\vartheta^3(13;25;46)(\W)\over \vartheta^3(15;23;46)(\W)},\\
\label{lambda6}
\ell_3&=&{\vartheta^3(13;26;45)(\W)\over \vartheta^3(16;23;45)(\W)}.
\end{eqnarray}
\end{proposition}

\begin{proposition}
\label{thetarel}
For the period matrix $\W$ of $C(\l),$ 
linear and cubic relations among $\vartheta^3(ij;kl;mn)(\W)$ 
coincide with the defining equations of $Y\subset \Ps^{14}:$
\begin{equation}
\label{lintheta}
\vartheta^3(ij;kl;mn)(\W)-\vartheta^3(ik;jl;mn)(\W)+\vartheta^3(il;jk;mn)(\W)=0,
\end{equation}
\begin{eqnarray}
\nonumber
& &\vartheta^3(ij;kl;mn)(\W)\vartheta^3(ik;jn;lm)(\W)\vartheta^3(im;jl;kn)(\W)
\\
\label{cubtheta}
&=&\vartheta^3(ij;kn;lm)(\W)\vartheta^3(ik;jl;mn)(\W)\vartheta^3(im;jn;kl)(\W).
\end{eqnarray}
\end{proposition}

Propositions \ref{invper} and \ref{thetarel} imply the following.
\begin{theorem}
\label{thetamap}
Let ${\mathit\Theta}$ be the map from $ \B^3/\G(1-\w)$ to $Y$ defined by 
$$x\mapsto [\dots,y_{\la ij;kl;mn\ra},\dots]=
[\dots,\vartheta^3(ij;kl;mn)(\jmath(x)),\dots].$$
We have the following $S_6$-equivariant commutative diagram:
$$
\matrix{ \L & \overset\psi\lto & \B^3/\G(1-\w) \cr
            &                  &               \cr
  \quad \leftset\iota\downarrow& \leftset{\mathit\Theta}\swarrow & \cr
            &                  &               \cr  
 \qquad Y\subset \Ps^{14}.  &    & \cr}
$$
\end{theorem}

In order to prove Propositions \ref{invper}, \ref{thetarel}, 
we state two facts in \cite{Mum};
the one is Riemann's theorem and the other is Abel's theorem.
\begin{fact}
\label{Riemann}
We suppose $z$ is a fix point on the Jacobi variety $Jac(R)$ 
of a Riemann surface $R$ of genus $r.$
The multi-valued function 
$\vartheta(z+\int_{P_0}^P\f,\tau)$ 
of $P$ on $X$ has $r$ zeros $P_1,\dots,P_r$ 
provided not to be constantly zero, 
where $\f=(\f_1,\dots,\f_r)$ is the normalized basis of the vector space of 
holomorphic $1$-forms on $R$ such that 
$(\int_{B_i}\f_j)_{ij}=I_r$
for a symplectic basis $\{A_1,\dots,A_r,B_1,\dots,B_r\}$ of $H_1(R,\Z),$ 
and $\tau=(\int_{A_i}\f_j)_{ij}.$
Moreover, there exists a point $\Delta$ on $Jac(R)$ called Riemann's constant 
such that 
$$z=\Delta-\sum_{i=1}^r \int_{P_0}^{P_i} \f.$$
\end{fact}
\begin{fact}
\label{abel}
Let $R$ be a Riemann surface of genus $r$ with an initial point $P_0.$ 
Suppose $\sum_{i=1}^d P_i$ and 
$\sum_{i=1}^d Q_i$ be effective divisors of degree $d$ satisfying 
\begin{equation}
\label{syuuki}
\sum_{i=1}^d\int_{P_0}^{P_i} \f = \sum_{i=1}^d\int_{P_0}^{Q_i} \f,
\end{equation}
where 
$\f$ is the normalized basis of vector space of holomorphic $1$-forms on $R.$ 
Then there exists a meromorphic function $f$ on $R$ such that 
$$(f)=\sum_{i=1}^d Q_i-\sum_{i=1}^d P_i;$$ 
$f$ can be expressed as 
$$f(P)=c{\prod_{i=1}^d\vartheta(e+\int_{Q_i}^P\f,\tau)\over 
\prod_{i=1}^d\vartheta(e+\int_{P_i}^P\f,\tau)},$$
where $c$ is a constant, $\tau$ is the period matrix of $R,$ 
$e$ satisfies $\vartheta(e)=0,$ 
$$\vartheta(e+\int_{P_i}^P\f,\tau)\not\equiv0,\quad 
\vartheta(e+\int_{Q_i}^P\f,\tau)\not\equiv0,$$
as multi-valued functions of $P$ on $R,$ 
and paths from $P_i$ and $Q_i$ to $P$ 
are the inverse of the paths in (\ref{syuuki})  followed by a common 
path from  $P_0$ to $P.$ 
\end{fact}
{\itshape Proof of Proposition \ref{invper}.} 
We take $R$ as 
$$C: w^3=z(z-1)(z-\ell_1)(z-\ell_2)(z-\ell_3)$$ 
with the initial point $P_0=(0,0)$ 
and put 
$$P_\infty=(\infty,\infty),\quad P_1=(1,0),\quad 
P_{\ell_i}=(\ell_i,0)\ (i=1,2,3).$$
Let us define a meromorphic function $f$ on $C$ by 
$C\ni (z,w)\mapsto z,$ then 
$$(f) = 3P_0-3P_\infty. $$

We construct a meromorphic function on $C$ with poles $3P_\infty$ 
and zeros $3P_0$ by following the recipe given in Fact \ref{abel}.
Let $\g_i(z_1,z_2)$ $(i=1,2,3)$ be a path in $C$ 
from $(z_1,w_1)$ to $(z_2,w_2)$ in the $i$-th sheet. 
Since $\w^2+\w+1=0,$ we have 
$$\sum_{i=1}^3 \int_{\g_i(0,\infty)}\f =(0,0,0,0)=
3\int_{P_0}^{P_0}\f$$
for three paths $\g_i(0,\infty)$ from $P_0$ to $P_\infty.$ 
We give the following table:
\begin{eqnarray*}
\int_{\g_1(\infty,0)}\f&=& {1\over 3}\int_{A_1-B_1}\f,\quad 
\int_{\g_2(\infty,0)}\f= {1\over 3}\int_{-2A_1-B_1}\f,\\
\int_{\g_3(\infty,0)}\f&=& {1\over 3}\int_{A_1+2B_1}\f,\quad 
\int_{\g_1(0,1)}\f= {1\over 3}\int_{-2A_1+A_2-A_4-B_1+2B_2+2B_4}\f,\\
\int_{\g_2(0,1)}\f&=& {1\over 3}\int_{A_1+A_2-A_4+2B_1-B_2-B_4}\f,\\
\int_{\g_3(0,1)}\f&=& {1\over 3}\int_{A_1-2A_2+2A_4-B_1-B_2-B_4}\f,\quad
\int_{\g_1(1,\ell_1)}\f= {1\over 3}\int_{A_2-B_2}\f,\\
\int_{\g_1(\ell_1,\ell_2)}\f&=& {1\over 3}\int_{-2A_2+A_3+2A_4-B_2+2B_3-B_4}\f,\quad
\int_{\g_1(\ell_2,\ell_3)}\f= {1\over 3}\int_{A_3-B_3}\f.
\end{eqnarray*}
Put 
$$e={1\over 6}\int_{3A_1+A_2+3A_3+5A_4-3B_1-B_2-3B_3+5B_4}\f,$$ 
corresponding to the characteristic ${1\over 6}(3,1,3,5)$ 
with label $(123),$
and define a meromorphic function $F$ of $P=(z,w)$ on $C$ as
\begin{equation}
\label{F}
F(P)={
\vartheta\left(e+\int_{\g_1(0,z)}\f,\W\right)^3\over
\prod_{i=1}^3\vartheta\left(e+\int_{\g_i(\infty,0)+\g_1(0,z)}\f,\W\right)},
\end{equation}
where $\W$ is the period matrix of $C.$
Since $\vartheta(123)$ vanishes, we have $\vartheta(e)=0.$ 
We check that neither the denominator nor the numerator of $F$ 
identically vanishes. We put  $P=P_{\ell_1},P_{\ell_2},P_{\ell_3}$ 
and use (\ref{thetaconst}) and (\ref{Z-shift}), then we have
\comment{
The $F(P)$ reduces to
$${\vartheta(a(\W-H),\W)^3\over \vartheta((a+{1\over3}e_1)(\W-H))^2
\vartheta((a-{2\over3}e_1)(\W-H))},
$$
where $a\in {1\over 6}\Z^4$ and $e_1=(1,0,0,0).$
Its denominator is 
$$\exp[-3\pi\sqrt{-1}(a\W\tr a-2aH\tr a)]\vartheta_{a,-aH}(\W)^3$$
and its numerator 
\begin{eqnarray*}
& &\exp[-2\pi\sqrt{-1}((a+{1\over3}e_1)(\W-2H)\tr(a+{1\over3}e_1))]
\vartheta_{a+{1\over3}e_1,-(a+{1\over3}e_1)H}(\W)^2\\
& &\cdot\exp[-\pi\sqrt{-1}((a-{2\over3}e_1)(\W-2H)\tr(a-{2\over3}e_1))]
\vartheta_{a-{2\over3}e_1,-(a-{2\over3}e_1)H}(\W)\\
&=&\exp[-3\pi\sqrt{-1}a(\W-2H)\tr a-{2\pi\sqrt{-1}\over 3}(\W_{11}-2)]\\
& &\cdot\vartheta_{a+{1\over3}e_1,-(a+{1\over3}e_1)H}(\W)^2
\vartheta_{a-{2\over3}e_1,-(a-{2\over3}e_1)H}(\W).
\end{eqnarray*}
Thus we have 
$$F(P)=\exp[{2\pi\sqrt{-1}/3}\W_{11}-{4\pi\sqrt{-1}/3}]
{\vartheta_{a,-aH}(\W)^3\over 
\vartheta_{a+{1\over3}e_1,-(a+{1\over3}e_1)H}(\W)^2
\vartheta_{a-{2\over3}e_1,-(a-{2\over3}e_1)H}(\W)}.$$}
\begin{eqnarray*}
F(P_{\ell_1})=cf(P_{\ell_1})=c\ell_1&=&\exp[{\pi\sqrt{-1}\over 3}(2\W_{11}+1)]
{\vartheta^3_{[-1,-1,3,-3]}(\W)\over \vartheta^3_{[1,-1,3,-3]}(\W)},\\
F(P_{\ell_2})=cf(P_{\ell_2})=c\ell_2&=&\exp[{\pi\sqrt{-1}\over 3}(2\W_{11}+1)]
{\vartheta^3_{[-1,1,-1,1]}(\W)\over \vartheta^3_{[1,1,-1,1]}(\W)},\\
F(P_{\ell_3})=cf(P_{\ell_3})=c\ell_3&=&\exp[{\pi\sqrt{-1}\over 3}(2\W_{11}+1)]
{\vartheta^3_{[-1,1,1,1]}(\W)\over \vartheta^3_{[1,1,1,1]}(\W)}, 
\end{eqnarray*}
where $c$ is a constant depending on $\W.$
By Proposition \ref{survive}, 
neither the denominator nor the numerator of $F$ identically vanishes. 

We put $P=P_{\infty},P_{0},P_{1};$ 
the denominator and the numerator of $F$ vanish at these points 
by Proposition \ref{survive}. 
Since $(F)=3P_0-3P_\infty,$ 
$P_\infty$ and $P_0$ are zeros of higher order of the denominator and 
numerator of $F,$ respectively. 
The number of zeros of the denominator and numerator of $F$ are $4$
by Fact ref{Riemann}, 
thus $P_1$ is a simple zero. 
We consider $\lim_{P\to P_1}F(P).$
Let $t$ be a local coordinate for $P$ around $P_1$ and 
$z(t)$ be $\int_{P_1}^P\f.$ 
We have 
$$F(P)=\exp[{\pi\sqrt{-1}\over 3}(2\W_{11}-2)]
{\vartheta^3_{[-1,-3,-3,-3]}(z(t),\W)\over \vartheta^3_{[1,3,3,3]}(z(t),\W)}.$$
When $P\to P_1,$ we have $t\to 0$ and $z(t)\to (0,0,0,0).$
Since $t=0$ is simple zero, we have 
$$
\lim_{t\to 0} 
{\vartheta^3_{[-1,-3,-3,-3]}(z(t),\W)\over \vartheta^3_{[1,3,3,3]}(z(t),\W)}
=
\lim_{t\to 0} 
{\vartheta^3_{[1,3,3,3]}(-z(t),\W)\over \vartheta^3_{[1,3,3,3]}(z(t),\W)}
=-1,$$
which implies $c=\exp[{\pi\sqrt{-1}\over 3}(2\W_{11}+1)].$
Hence we have the expressions (\ref{lambda4}), (\ref{lambda5}) and 
(\ref{lambda6}). \qed

{\itshape Proof of Proposition \ref{thetarel}.} 
In order to obtain a cubic relation among $\vartheta^3(ij;kl;mn)$'s, 
put 
$$e={1\over 6}\int_{3A_1+5A_2+3A_3+5A_4-3B_1-5B_2-3B_3+5B_4}\f,$$ 
corresponding to the characteristic ${1\over 6}(3,5,3,5)$ 
with label $(124),$ then $\vartheta(e)=0;$ 
and define a meromorphic function $F$ by (\ref{F}).
We have 
\begin{eqnarray*}
F(P_1)=cf(P_1)=c&=&\exp[{\pi\sqrt{-1}\over 3}(2\W_{11}+1)]
{\vartheta^3_{[-1,1,3,-3]}(\W)
\over \vartheta_{[1,1,3,-3]}(\W)},\\
F(P_{\ell_2})=cf(P_{\ell_2})=c\ell_2&=&\exp[{\pi\sqrt{-1}\over 3}(2\W_{11}+1)]
{\vartheta^3_{[-1,-1,-1,1]}(\W)\over \vartheta^3_{[1,-1,-1,1]}(\W)},\\
F(P_{\ell_3})=cf(P_{\ell_3})=c\ell_3&=&\exp[{\pi\sqrt{-1}\over 3}(2\W_{11}+1)]
{\vartheta^3_{[-1,-1,1,1]}(\W)\over \vartheta^3_{[1,-1,1,1]}(\W)},
\end{eqnarray*}
and 
\begin{eqnarray*}
& &c\ell_1=cf(P_{\ell_1})=\lim_{P\to P_{\ell_1}}F(P)\\
&=&\exp[{\pi\sqrt{-1}\over 3}(2\W_{11}-2)]
{\vartheta^3_{[-1,-3,-3,-3]}(\int_{P_{\ell_1}}^P \f,\W)\over 
\vartheta^3_{[1,3,3,3]}(\int_{P_{\ell_1}}^P\f,\W)}=
\exp[{\pi\sqrt{-1}\over 3}(2\W_{11}+1)].\\
\end{eqnarray*}
These imply 
\begin{eqnarray*}
\ell_1={cf(P_{\ell_1})\over cf(P_1)}&=&{\vartheta^3(13;24;56)(\W)\over 
\vartheta^3(14;23;56)(\W)},\\
\ell_2={cf(P_{\ell_2})\over cf(P_1)}&=&
{\vartheta^3(14;25;36)(\W)\vartheta^3(13;24;56)(\W)
\over \vartheta^3(15;24;36)(\W)\vartheta^3(14;23;56)(\W)},\\
\ell_3={cf(P_{\ell_3})\over cf(P_1)}
&=&{\vartheta^3(14;26;35)(\W)\vartheta^3(13;24;56)(\W)
\over \vartheta^3(16;24;35)(\W)\vartheta^3(14;23;56)(\W)}.
\end{eqnarray*}
Compare with the above expression of $\ell_2$ and (\ref{lambda5}), 
we have a cubic relation among the $\vartheta^3(ij;kl;mn)$'s.
By letting $\tilde\G$ act on theta constants, we have more cubic relations  
among $\vartheta^3(ij;kl;mn)$'s. 

Let us lead a linear relation among the $\vartheta^3(ij;kl;mn)$'s. 
We start with the meromorphic function $f':(z,w)\mapsto z-1;$
note that $(f')=3P_1-3P_\infty.$ 
Put 
$$e={1\over 6}\int_{3A_1+A_2+3A_3+5A_4-3B_1-B_2-3B_3+5B_4}\f,$$ 
corresponding to the characteristic ${1\over 6}(3,1,3,5)$ 
with label $(123),$
and define a meromorphic function $F'$ of $P=(z,w)$ on $C$ as
$$F'(P)={\prod_{i=1}^3
\vartheta\left(e+\int_{\g_i(1,0)+\g_1(0,z)}\f,\W\right)\over 
\prod_{i=1}^3
\vartheta\left(e+\int_{\g_i(\infty,0)+\g_1(0,z)}\f,\W\right)
}.$$
Since $\vartheta(123)$ vanishes, we have $\vartheta(e)=0.$ 
We consider $\lim_{P\to P_0}F'(P)$ and 
put $P=P_{\ell_1}$ then we have
\begin{eqnarray*}
F'(P_0)&=&cf'(P_0)=-c=\lim_{P\to P_0}K\exp[{4\pi\sqrt{-1}\over 3}]
{\vartheta^3_{[-5,-1,-3,-5]}(\int_{P_0}^P\f,\W)\over 
\vartheta^3_{[5,1,3,5]}(\int_{P_0}^P\f,\W)}\\
F'(P_{\ell_1})&=&cf'(P_{\ell_1})=c(\ell_1-1)=K\exp[{4\pi\sqrt{-1}\over 3}]
{\vartheta^3_{[3,3,3,-1]}(\W)\over \vartheta^3_{[1,-1,3,-3)}(\W)},
\end{eqnarray*}
where 
$$K=\exp[-{2\pi\sqrt{-1}\over 3}e'\W\tr(e'-e_1)+
{4\pi\sqrt{-1}\over3}e'H\tr(e'-e_1)]$$
and $e'=(1,-1,0,1).$ Now we have the expression 
$$\ell_1-1={\vartheta^3(12;34;56)(\W)\over \vartheta^3(14;23;56)(\W)}.$$
Since we had in (\ref{lambda4})
$$\ell_1={\vartheta^3(13;24;56)(\W)\over \vartheta^3(14;23;56)(\W)},$$
we get a relation 
$${\vartheta^3(12;34;56)(\W)\over \vartheta^3(14;23;56)(\W)}
-{\vartheta^3(13;24;56)(\W)\over \vartheta^3(14;23;56)(\W)}+1=0,$$
which is equivalent to 
$$\vartheta^3(12;34;56)(\W)-\vartheta^3(13;24;56)(\W)
+\vartheta^3(14;23;56)(\W)=0.$$
Action of $\tilde\G$ produces other the linear relations
among the $\vartheta^3(ij;kl;mn)$'s.  \qed

\section{Appendix}
In this section, we give a geometrical meaning of the label of $a$'s.
In order to do this, we determine  Riemann's constant $\Delta.$ 
\begin{fact}
Riemann's constant $\Delta$ is given  by 
\begin{equation}
\label{R-const}
\Delta=\sum_{i=1}^{m+r-1} \int_{P_0}^{P_i}\f-\sum_{j=1}^m\int_{P_0}^{Q_j}\f
\end{equation}
for a certain divisor $D_0=\sum_{i=1}^{m+r-1}P_i-\sum_{j=1}^{m}Q_j$ 
such that $2D_0$ is linearly equivalent to the canonical divisor of $R.$ 
It is easy to see that Riemann's constant $\Delta$ is 
a half period on $Jac(R)$ if and only if $(2r-2)P_0$  is a canonical divisor.
\end{fact}

For our case, Riemann's constant $\Delta$ is a half period on $Jac(C(\l))$ 
since we have $6P_0=(\f_4)$ for any $C(\l).$  

\begin{proposition}
Riemann's constant $\Delta$ is invariant under the action of 
the monodromy group $\tilde\G(1-\w).$ Hence we have  
$$\Delta=({1\over 2},\dots,{1\over 2}).$$
\end{proposition}
{\itshape Proof.}
Let $\g$ be a closed path in $\L$ and $g\in \tilde\G(1-\w)$ be 
its representation. Since $\Delta$ is a half period point of $Jac(C(\l)),$ 
it is expressed by $c=(c_1,\dots,c_8)$ $(c_i\in \{0,1/2\}).$ 
When $\l$ moves a little, this vector is invariant and  presents $\Delta.$
By the continuation along $\g,$ $\Delta$ is presented by the vector $c$ 
with respect to the transformed homology basis by $g;$ 
i.e., it is presented by $g\cdot c$ with respect to the initial homology 
basis. 

On the other hand, $\Delta$ is invariant as a point of $Jac(C(\l))$ 
under the continuation along $\g$ with respect to the initial basis 
by the expression (\ref{R-const}). 
Thus we have $g\cdot c=c.$ There is only one half characteristic  
$({1\over 2},\dots,{1\over 2})$ invariant under $\tilde\G(1-\w).$ 
\qed

By straightforward calculation, we have the following proposition 
giving a geometrical meaning of the label of $a$'s.
\begin{proposition}
The points $(a,-aH)$ of $Jac(C)$ for $a$ with label $(ijk)$ and $(i^2j)$  are 
expressed as 
\begin{eqnarray*}
& &\Delta-\int_{P_0}^{P_{\l_i}}\f-\int_{P_0}^{P_{\l_j}}\f-
\int_{P_0}^{P_{\l_k}}\f,\\
& &\Delta-2\int_{P_0}^{P_{\l_i}}\f-\int_{P_0}^{P_{\l_j}}\f,
\end{eqnarray*}
respectively.
\end{proposition}

We have the necessary and sufficient condition for $\vartheta(z,\tau)=0.$ 

\begin{fact}
\label{zerotheta}
For a period matrix $\tau$ of Riemann's surface $R$ of genus $r,$ 
$\vartheta(z,\tau)=0$ if and only if there exists an effective divisor 
$\sum_{i=1}^{r-1}P_i$ such that 
$$z=\Delta-\sum_{i=1}^{r-1}\int_{P_0}^{P_i}\f.$$
\end{fact}

\begin{proposition}
\label{zeros}
The theta constant $\vartheta(ij;kl;mn)(\jmath(x))$ 
vanishes only on the $\G(1-\w)$ orbit of 
the mirrors of $R^\w_{ij},$ $R^\w_{kl}$ and $R^\w_{mn}.$
\end{proposition}
{\itshape Proof.} The function $\vartheta(13;24;56)(\jmath(x))$  is a non-zero 
constant times 
$$\vartheta(\Delta-\int_{P_0}^{P_\infty}\f-\int_{P_0}^{P_0}\f-
\int_{P_0}^{P_1}\f+\int_{P_0}^{P_{\ell_1}}\f,\jmath(x)).$$
By the previous fact, $\vartheta(13;24;56)(\jmath(x))=0$ if and only if  
there exists an effective divisor $Q_1+Q_2+Q_3$ such that 
$$Q_1+Q_2+Q_3\equiv P_\infty+2P_0+P_1-P_{\ell_1}=E.$$ 
By the Riemann-Roch theorem, 
the dimension of vector space of  meromorphic functions $f$
such that $(f)+E\ge 0$ is equal to 
that of meromorphic $1$-forms $\phi$ such that 
\begin{equation}
\label{Roch}
(\phi) -E\ge 0.
\end{equation}
Since we have
$$(\f_1)=P_\infty+P_0+P_1+P_{\ell_1}+P_{\ell_2}+P_{\ell_3},\quad
(\f_2)=6P_\infty,$$
$$(\f_3)=3P_\infty+3P_0,\quad (\f_4)=6P_0,$$
there does not exist a meromorphic $1$-from satisfying (\ref{Roch}).
Thus if $\l\in \L$ then 
no effective divisor $Q_1+Q_2+Q_3$ such that $Q_1+Q_2+Q_3\equiv E.$ 

The zeros of theta constants on mirrors are studied in \cite{Shi},  
which yields this proposition.
\qed


\end{document}